\documentclass[psamsfonts]{amsproc}
\RequirePackage{amssymb}
\RequirePackage{epsfig,graphics}
\RequirePackage{upref}
\newtheorem{theorem}{Theorem}[section]
\newtheorem{lemma}[theorem]{Lemma}
\newtheorem{proposition}[theorem]{Proposition}
\newtheorem{corollary}[theorem]{Corollary}
\theoremstyle{definition}
\newtheorem{definition}[theorem]{Definition}

\theoremstyle{remark}
\newtheorem{remark}[theorem]{Remark}
\newtheorem{algorithm}{Algorithm}
\numberwithin{equation}{section}

\newcommand{\RR}{\mathbb{R}}
\newcommand{\NN}{\mathbb{N}}

\begin{document}

\title{A bisection algorithm for the numerical Mountain Pass}
\author[V. Barutello]{Vivina Barutello}
\author[S.Terracini]{Susanna Terracini}
\address{Dipartimento di Matematica e Applicazioni, Universit\`{a} di Milano-Bicocca,
         via Cozzi 53, 20125 Milano, Italy.}
\email{vivina@matapp.unimib.it, suster@matapp.unimib.it}
\thanks{This work is partially supported by M.I.U.R. project 
      ``Metodi Variazionali ed Equazioni Differenziali Nonlineari".}

\begin{abstract}
We propose a constructive proof for the Ambrosetti-Rabinowitz Mountain Pass Theorem
providing an algorithm, based on a bisection method, for its implementation.
The efficiency of our algorithm, particularly suitable for problems in high dimensions,
consists in the low number of flow lines to be computed for its convergence;
for this reason it improves the one currently used and proposed by Y.S. Choi and 
P.J. McKenna in~\cite{CMc}.
\end{abstract}

\maketitle

\section*{Introduction}
\label{sec:intro}

This paper deals with constructive methods to seek critical points 
which are not minimizers for functionals defined on Hilbert spaces. 
The existence of such critical points may be detected from the topological features 
of the sublevels, under some compactness conditions.
The purpose of this work is to give a constructive version of the 
Ambrosetti-Rabinowitz Mountain Pass Theorem (see~\cite{AR}), which is one of the most useful 
abstract critical point theorems and that has found relevant 
applications in solving nonlinear boundary value problems.
Our constructive proof relies on an algorithm based on a bisection method;
this algorithm, elaborated in different versions depending on the topology
of the sublevels (see Algorithms~\ref{alg:1b} and~\ref{alg:2}), 
can be implemented numerically.
In spite of its simplicity, the algorithm proposed is new and in terms of numerical optimization,
it improves the one currently used and proposed by Y.S. Choi and P.J. McKenna in~\cite{CMc} 
(see also~\cite{CMcR}).
The efficiency of our algorithm consists in the low number of steepest descent 
flow lines to be determined to obtain a good approximation of a critical point of mountain pass type; 
this fact makes the method particularly fit for problems in high dimensions such as those coming from 
the discretization of nonlinear boundary value problems in infinite dimensional spaces.
Moreover, compared with the usual Newton's method, our algorithm has two main advantages:
it always converges and it does not need an a priori good initial guess.

For some numerical applications of Algorithms~\ref{alg:1b} and~\ref{alg:2} we refer the reader to Chapter 2 of~\cite{PhD_vi}, where such algorithms are applied to the action functional associated to
the $n$-body problem with simple or double choreography constraint respectively.
The theory exposed in this paper, permits the determination of a new solution for the 
$3$-body problem in a rotating frame with angular velocity $\omega = 1.5$.

\section{An iterative algorithm for critical points}
\label{sec:1}
\subsection{The steepest descent flow}
\label{subsec:steep}
We consider a Hilbert space $X$ and a functional $f:X \rightarrow \RR$ of class $C^2$. 
Fixed $c \in \RR$, we define the $c$-\emph{sublevel} of $f$ as
the open set
$$
f^c := \{ x \in X : f(x) < c \}
$$ 
and the set of critical points of $f$ as
$$
Crit(f):=\{ x \in X : \nabla f(x) = 0 \}.
$$
The point $x_0 \in Crit(f)$ is a \emph{local minimizer} 
for the functional $f$ if there exists $r>0$, such that 
$f(x) \geq f(x_0)$, $\forall x \in B_r(x_0)$;
$x_0$ is a \emph{strict local minimizer} if there exists $r_0>0$ such that 
for every $r < r_0$, $\displaystyle \inf_{x \in \partial B_r(x_0)}f(x) > f(x_0)$.

Let $\eta:\RR_+ \times X \rightarrow X$ be the steepest descent flow 
associated with the functional $f$ defined as the solution of the Cauchy problem
\begin{equation}
\label{flow}
\left\{ 
\begin{array}{l}
\displaystyle\frac{d}{dt}\eta(t,x) = -\frac{\nabla f(\eta(t,x))}{1+\|\nabla f(\eta(t,x))\|} \\ \eta(0,x) = x
\end{array}
\right.
\end{equation}
We say that a subset $X_0 \subset X$ is \emph{positively invariant} 
for the flow $\eta$ if $\{ \eta(t,x_0), t \geq 0 \} \subset  X_0$, for every $x_0 \in X_0$.
We term $\omega$-\emph{limit of} $x \in X$ for the flow $\eta$, 
the closed positively invariant set
$$
\omega_{x} = \left\{ \lim_{t_n\rightarrow +\infty}\eta(t_n,x): (t_n)_n \subset \RR_+ \right\}
$$
We now state two useful preliminary Lemmata concerning
some properties of the steepest descent flow defined in (\ref{flow}); 
for their simple proofs we refer, if necessary, to~\cite{PhD_vi}.

\begin{lemma}
\label{mis_tempi}
Let $\eta$ be the steepest descent flow defined in~(\ref{flow}); let $x \in X$,
$\gamma \in (0,1]$ and $T>0$, then
$$
\left| \left\{ t \in [0,T]:\| \nabla f(\eta(t,x)) \|\geq \gamma \right\} \right| 
\leq \frac{f(x) - f(\eta(T,x))}{\gamma^2 / 2}.
$$
\end{lemma}

\begin{lemma}
\label{omega_limit}
Let $\eta$ be the steepest descent flow defined in~(\ref{flow}); let $x \in X$,
then
\begin{itemize}
\item[(i)] $\omega_x \subset Crit(f)$;
\end{itemize}
moreover if $ c_x := \lim_{t\rightarrow +\infty} f(\eta(t,x)) > -\infty$, then
\begin{itemize}
\item[(ii)]  
$\displaystyle \exists (t_n)_n, t_n \rightarrow +\infty$ such that 
$$
\lim_{n \rightarrow +\infty}f(\eta(t_n,x))=c_x \mbox{ and  }
\displaystyle\lim_{n \rightarrow +\infty} \nabla f(\eta(t_n,x))=0.
$$
\end{itemize}
\end{lemma}

\subsection{Disconnected sublevels}
\label{subsec:disc}

Let $c \in \RR$ be such that the sublevel $f^c$ is disconnected,
we term $(F^{c}_i)_i$ its disjoint connected components
$$
f^{c} = \bigcup_{i} F^{c}_i, \quad F^{c}_i \cap F^{c}_j = \emptyset, \,\, \forall i \neq j.
$$ 
For every index $i$, we consider the basin of attraction 
of the set $F^{c}_i$
$$
\mathcal{F}^{c}_i := \left\{ x \in X : \omega_x \subset F^{c}_i \right\}.
$$

\begin{proposition} 
\label{prop cal F}
Let $f^c$ be a disconnected sublevel for the functional $f$.
Let $F_i^c$ be the disjoint connected components of $f^c$ and 
$\mathcal{F}^{c}_i$ their basins of attraction. Then, for every  index $i$,
the following assertions hold:
\begin{enumerate}
\item[(i)] $\mathcal{F}^c_i$ is an open set;
\item[(ii)] $\partial \mathcal{F}^c_i$ is a positively invariant set;
\item[(iii)] $\displaystyle \inf_{x \in \partial \mathcal{F}^c_i} f(x) \geq c$;   
\item[(iv)] minimizers of $f$ in $\partial \mathcal{F}^c_i$  are critical points for $f$;
\item[(v)] critical points for $f$ in the set $\partial\mathcal{F}^c_i$ are not strict local minimizers.
\end{enumerate}
\end{proposition}
\begin{proof}
We prove properties {(i)}-{(v)} when $i=1$.

{(i)} Let $x \in \mathcal{F}^c_1$, then $\omega_x \subset F^c_1$
and, since $\eta$ is a gradient flow and $F^c_1$ is a connected component 
of the sublevel $f^c$, there exists $T$ such that $\eta(t,x) \in F^c_1$, 
for every $t \geq T$. The continuity of the flow $\eta$ ensures that there exists
$\delta>0$ such that $\eta(T,B_\delta(x))\subset  F^c_1$ and, hence, 
$B_\delta(x) \subset \mathcal{F}^c_1$. 

{(ii)} By the sake of contradiction, suppose there exist 
$x \in \partial\mathcal{F}^c_1$ and $T>0$ such that 
$x_T := \eta(T,x) \notin \partial \mathcal{F}^c_1$; if $x_T \in \mathcal{F}^c_1$,
then there exists $\bar T$ such that $\eta(\bar T,x) \in F^c_1$, hence 
$x \in \mathcal{F}^c_1$, but this contradicts {(i)}.
If $x_T \in X \backslash\bar{\mathcal{F}}^c_1$, then there exists $\epsilon >0$ such that 
$B_\epsilon(x_T) \subset X \backslash \bar {\mathcal{F}}^c_1$;
hence there exists $\delta >0$ such that $\eta(T,B_\delta(x)) \subset B_\epsilon(x_T)$
in contradiction with the definition of the set $\mathcal{F}^c_1$.

{(iii)} Follows from the definition of the set $\mathcal{F}^c_1$.
 
{(iv)} Follows from {(ii)} and $\omega_x \subset Crit(f)$, for every $x \in \partial\mathcal{F}^c_i$.

{(v)} Let $\bar x \in \partial\mathcal{F}^c_1$ be a critical point for 
the functional $f$; by the sake of contradiction, suppose that $\bar x$ is a strict local minimizer,
hence $\bar x \notin \partial{F}^c_1$.
Let $r>0$ be such that $B_r(\bar x) \cap F_1^c = \emptyset$
and $\alpha_r>0$ such that 
$$
f(x) \geq \alpha_r+f(\bar x), \quad \forall x \in \partial B_r(\bar x).
$$
By the continuity of $f$, the set
$$
U := \{ x \in B_r(\bar x) : f(x)<\alpha_r+f(\bar x) \},
$$ 
is a neighborhood of $\bar x$, hence there exists a point
$x_U \in U \cap \mathcal{F}_1^c$; since the flow $\eta$ is dissipative, 
and $f(\partial B_r(\bar x)) > f(x_U)$, we claim that
$\omega_{x_U} \subset B_r(\bar x)$, which contradicts the definition of $\mathcal{F}_1^c$.
\end{proof}

We term \emph{path} a continuous function $\gamma : [0,1] \rightarrow X$. 
Given a pair of points $x_1,x_2 \in X$, $x_1 \neq x_2$, we define the set 
of paths joining $x_1$ to $x_2$ as 
\begin{equation}
\label{Gamma_set}
\Gamma_{x_1,x_2} := \left\{ \gamma \in C([0,1],X) : \gamma(0)=x_1, \gamma(1)=x_2 \right\}.
\end{equation}

\begin{theorem}
\label{thm:alg1}
Let $f^c$ be a disconnected sublevel for the functional $f$.
Let $F_i^c$ be the disjoint connected components of $f^c$ and 
$\mathcal{F}^{c}_i$ their basins of attraction. 
Let $x_i \in F_i^c$, $i=1,2$, and $\gamma \in \Gamma_{x_1,x_2}$;
then there exists $\bar x \in \gamma([0,1]) \cap \partial \mathcal{F}^{c}_1$. 
\end{theorem}

\begin{proof}
The first step is the description of an algorithm
that, given a path $\gamma$ in the set  $\Gamma_{x_1,x_2}$, selects a point 
$\bar x \in \gamma([0,1]) \cap \partial \mathcal{F}^{c}_1$.

\begin{algorithm}
\label{alg:1a}
\texttt{
\begin{itemize}
\item[{ Step 0.}] $s_1^0 = 0$, $s_2^0 = 1$, $\displaystyle s_m^0=\frac{s^0_1+s^0_2}{2}$, 
                  $x_1^0 = x_1$, $x_2^0 = x_2$, $x_m^0 = \gamma (s_m^0)$
\item[{ Step i.}] if $\omega_{x_m^{i-1}} \subset F^{c_0}_1$,
									\hspace{.5cm} $s_1^{i}=s_m^{i-1}$, $s_2^i=s_2^{i-1}$\\
                  else $s_1^i=s_1^{i-1}$, $s_2^i=s_m^{i-1}$\\
                  $x_1^i = \gamma(s_1^i)$, $x_2^i = \gamma(s_2^i)$,
                  $\displaystyle s_m^i = \frac{s^i_1+s^i_2}{2}$,
                  $x_m^i = \gamma (s_m^i)$
\end{itemize}}
\end{algorithm}
We have then defined two sequences $(s_1^i)_i,(s_2^i)_i$ such that
$$
0=s_1^0 \leq s_1^1 \leq \ldots \leq s_1^i < s_2^i \leq \ldots \leq s_2^1 \leq s_2^0 =1  
$$
and
$$
|s_1^i-s_2^i|=\frac{1}{2^i} \stackrel{i \rightarrow +\infty}{\longrightarrow} 0.
$$
Since $(s_1^i)_i,(s_2^i)_i$ are bounded monotone sequences
$$
\lim_{i \rightarrow + \infty}s_1^{i} =  \lim_{i \rightarrow + \infty} s_2^{i} = {\bar s} \in [0,1].
$$
Let ${\bar x}=\gamma({\bar s})$, hence 
$$
\lim_{i \rightarrow + \infty}x_1^{i} =  \lim_{i \rightarrow + \infty} x_2^{i} = {\bar x},
$$
and, necessarily ${\bar x}$ lies on the positively invariant set $\partial\mathcal{F}^{c}_1$. 
\end{proof}

\begin{corollary}
\label{cor:thm:alg1}
In the same conditions of Theorem~\ref{thm:alg1}, let
$\bar x \in \gamma([0,1]) \cap \partial \mathcal{F}^{c}_1$,
then $f(\omega_{\bar x}) \geq {c}$ and there exists a sequence
$(x_n)_n = \eta(t_n, \bar x) \subset \partial \mathcal{F}^{c}_1$ such that
$$
\lim_{n \rightarrow +\infty}\nabla f (x_n)=0,  \quad
\lim_{n \rightarrow +\infty}f (x_n)=f(\omega_{\bar x}).
$$
\end{corollary}
\begin{proof}
Since $\partial\mathcal{F}^{c}_1$ is a positively invariant closed set,  
$\omega_{\bar x} \subset \partial\mathcal{F}^{c}_1$ and $f(\omega_{\bar x}) \geq c$.
We use Lemma~\ref{omega_limit}, to define the sequence $(x_n)_n$.
\end{proof}

Unfortunately, the proof of Corollary~\ref{cor:thm:alg1} is not constructive 
in the sense that it does not provide a method to determine the critical 
set $\omega_{\bar x}$. The reasons why we can not have an implementable 
proof of this result are, first, that we can not determine precisely, 
in a finite number of steps, the point $\bar x$, since
it is the limit of the sequence $(x_1^n)_n$ in Algorithm~\ref{alg:1a}. Second, 
$\omega_{\bar x}$ is defined as a limit for $t \rightarrow +\infty$ and 
we are not able to determine the value $f(\omega_{\bar x})$ numerically.

Although the following result is a consequence of Corollary~\ref{cor:thm:alg1},
its relevance consists in its constructive proof, that gives a method 
to determine a critical point for the functional $f$ at a level higher than $c$.

\begin{corollary}
\label{cor2:thm:alg1}
In the same conditions of Theorem~\ref{thm:alg1}, let
$\bar x \in \gamma([0,1]) \cap \partial \mathcal{F}^{c}_1$,
then there exists a sequence $(\tilde y_n)_n\subset X$, 
$\tilde y_n:= \eta(\tilde T_n,x_1^n)$, such that
$$
\lim_{n \rightarrow +\infty}\nabla f (\tilde y_n)=0, \mbox{ and } 
c \leq f (\tilde y_n) \leq f({\bar x}),\, \forall n \in \NN. 
$$
\end{corollary}
\begin{proof}
Let $(x_1^n)_n$ be the sequence defined in Algorithm~\ref{alg:1a},
$(x_1^n)_n \subset \mathcal{F}_1^c \cap \gamma([0,1])$ 
and converging to $\bar x \in \partial \mathcal{F}_1^c$.
Since $(x_1^n)_n \subset \mathcal{F}_1^c$, we can define 
$$
T_n:= \inf \left\{ t \geq 0 : f(\eta(t,x_1^n)) \leq c \right\}
$$
and 
$$
\gamma_n:= \inf_{t \in [0,T_n]} \| \nabla f(\eta(t,x_1^n))\|.
$$
Hence, using Lemma~\ref{mis_tempi} (we can suppose $\gamma_n < 1$) we deduce
$$
T_n = \left| \left\{ t \in [0,T_n] : \| \nabla f(\eta(t,x_1^n))\|\geq \gamma_n \right\}\right|
    \leq \frac{2(f(x_1^n)-c)}{\gamma_n^2}
$$
and we can conclude that 
$$
\gamma_n \leq \sqrt{ \frac{2(f(x_1^n)-c)}{T_n}}.
$$
Since $(x_1^n)_n  \rightarrow \bar x$ as $n \rightarrow +\infty$
and $f(\omega_{\bar x}) \geq c$, we have
$$
\lim_{n \rightarrow +\infty}T_n = +\infty \,\, \mbox{ and } \lim_{n \rightarrow +\infty}\gamma_n = 0. 
$$
We deduce the existence of a sequence $(\tilde T_n)_n$ such that 
$(\tilde T_n)_n \in [0,T_n]$ and, defining $\tilde y_n := \eta(\tilde T_n,x_1^n)$,
we have 
$$
\lim_{n \rightarrow +\infty} \| \nabla f(\tilde y_n)\| = 0.
$$
\end{proof}

The proof of Corollary~\ref{cor2:thm:alg1} shows that Algorithm~\ref{alg:1a} can be 
improved in the following Algorithm~\ref{alg:1b} to obtain the sequence $(\tilde y_n)_n$.

\begin{algorithm}
\label{alg:1b}
\texttt{
\begin{itemize}
\item[{ Step 0.}] $s_1^0 = 0$, $s_2^0 = 1$, $\displaystyle s_m^0=\frac{s^0_1+s^0_2}{2}$, 
                  $x_1^0 = x_1$, $x_2^0 = x_2$, $x_m^0 = \gamma (s_m^0)$
\item[{ Step i.}] if $\omega_{x_m^{i-1}} \subset F^{c_0}_1$,
									\hspace{.5cm} $s_1^{i}=s_m^{i-1}$, $s_2^{i}=s_2^{i-1}$\\
                  else $s_1^{i}=s_1^{i-1}$, $s_2^{i}=s_m^{i-1}$\\
                  $x_1^{i} = \gamma(s_1^{i})$, $x_2^{i} = \gamma(s_2^{i})$,
                  $\displaystyle s_m^{i} = \frac{s^{i}_1+s^{i}_2}{2}$,
                  $x_m^{i} = \gamma (s_m^{i})$\\
                  $\displaystyle T_{i} := \inf \left\{ t \geq 0 : f(\eta(t,x_1^{i})) \leq c \right\}$\\
                  $\displaystyle \tilde T_{i} := t \in [0,T_{i}]$ that minimizes $\| \nabla 
                  f(\eta(t,x_1^{i}))\|$\\
                  $\displaystyle \tilde y_{i} := \eta(\tilde T_{i},x_1^{i}).$
\end{itemize}}
\end{algorithm}
The sequence $(\tilde y_n)_n$ defined in Algorithm~\ref{alg:1b} will converge, 
once we impose some additional compactness conditions on the functional $f$. 
In this sense we give the following definitions.

\begin{definition}
\label{def:PS}
A sequence $(x_m)_m \subset X$ is termed a 
\emph{Palais-Smale sequence in the interval $[a,b]$ for the functional $f$} if
$$
a \leq f(x_m) \leq b, \forall m \in \NN \quad \mbox{ and } \quad
\nabla f(x_m)\stackrel{m \rightarrow +\infty}{\longrightarrow} 0.
$$
The functional \emph{$f$ satisfies the Palais-Smale condition in the interval $[a,b]$}
if every Palais-Smale sequence in the interval $[a,b]$ for the functional $f$,
$(x_m)_m$, has a  converging subsequence $x_{m_k} \rightarrow x_0 \in X$.
Similarly, a sequence $(x_m)_m \subset X$ is a 
\emph{Palais-Smale sequence at level $c$ for the functional $f$} if
$$
f(x_m)\stackrel{m \rightarrow +\infty}{\longrightarrow} c \quad \mbox{ and } \quad
\nabla f(x_m)\stackrel{m \rightarrow +\infty}{\longrightarrow} 0.
$$
The functional \emph{$f$ satisfies the Palais-Smale condition at level $c$, 
(PS)$_{c}$}, if every Palais-Smale sequence at level $c$ for $f$
has a  converging subsequence.
\end{definition}

\begin{remark}
Corollary~\ref{cor:thm:alg1} ensures the existence of
a Palais-Smale  sequence at level $f(\omega_{\bar x})$ for the functional $f$. 
When the functional $f$  satisfies the (PS)$_{f(\omega_{\bar x})}$, 
we conclude that there exists a  critical point $\bar{\bar x}$ for $f$ 
such that $f(\bar{\bar x}) = f(\omega_{\bar x})$. 
Corollary~\ref{cor2:thm:alg1} implies the existence of
a Palais-Smale sequence  for $f$ in the interval $[c,f(\bar x)]$. 
When the functional verifies the Palais-Smale condition in $[c,f(\bar x)]$,
it ensures the convergence of the sequence $(\tilde y_n)_n$
constructed in Algorithm~\ref{alg:1b}.
\end{remark}

\begin{figure}[ht!]
\begin{center}
\psfig{figure=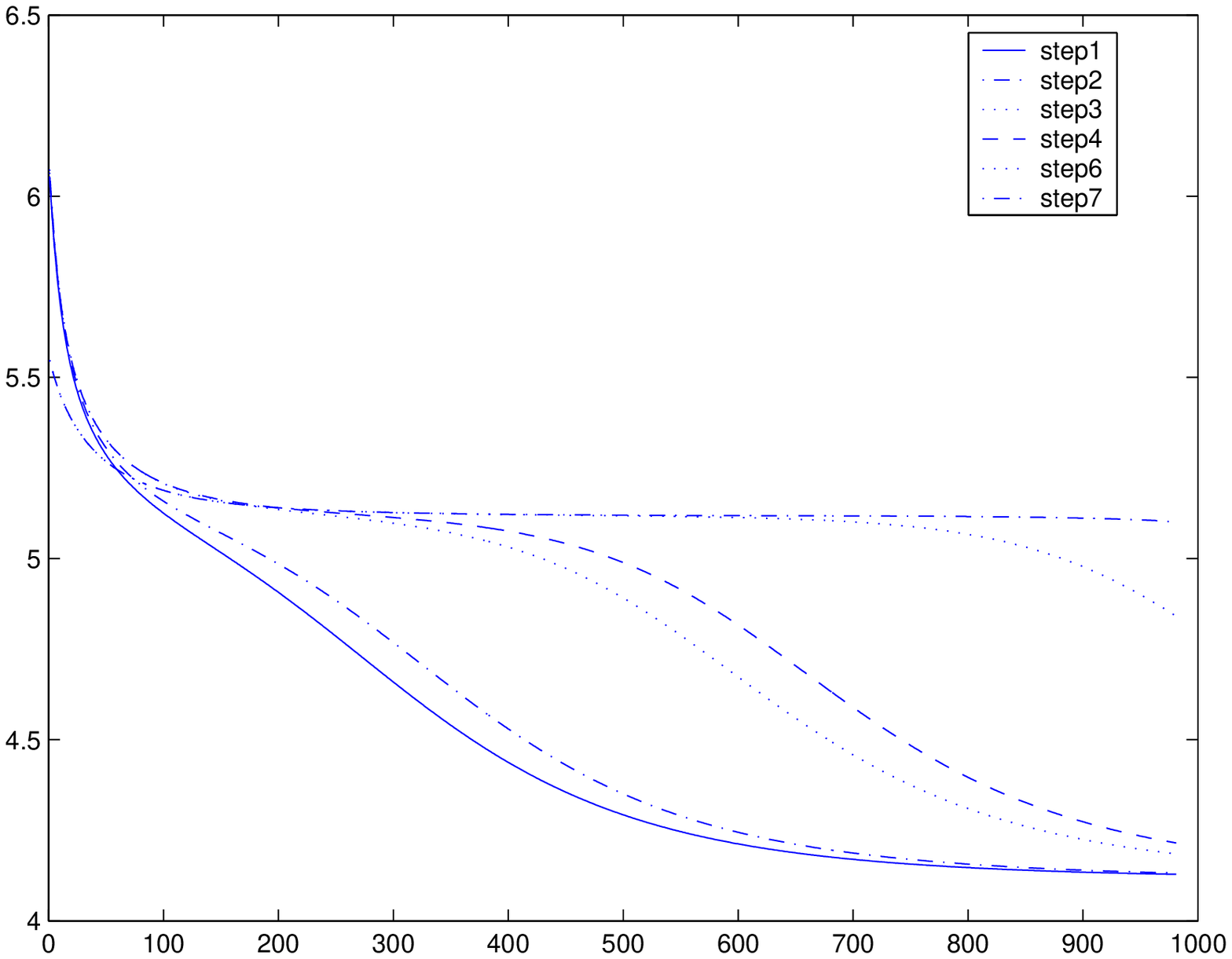,width=9.0cm}\\
\psfig{figure=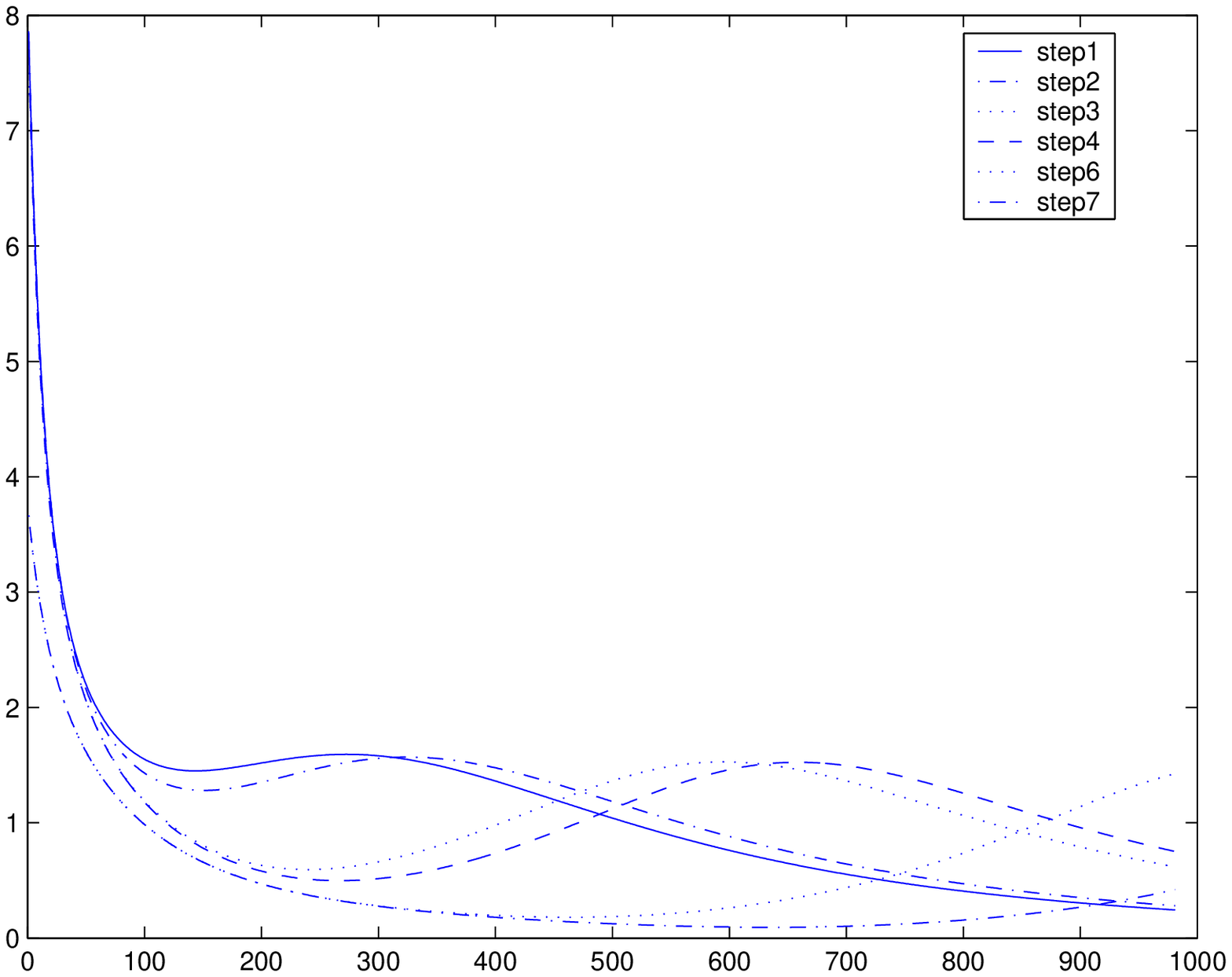,width=9.0cm}
\end{center}
\caption{
\label{fig:action_levels}
Action functional and norm of its gradient on the descent 
flow from the first elements of the sequence $(x_1^j)_j$.}
\end{figure}

Figure~\ref{fig:action_levels} shows the implementation of Algorithm~\ref{alg:1b}.
In the first graph, each line represents 
the values of the action functional on the steepest descent flow 
starting at the first elements  of the sequence $(x^i_1)_i$. 
Remark that, as the index $i$
increases, the interval of time in which the action function 
is approximately~$5.2$ increases, too. 
This is due to the fact that the line of the flow departing  
from $x^i_1$ passes each step closer to the desired critical point, whose action level 
is approximately~$5.2$.
In the second picture, we see the value of the norm of the gradient of
the action functional on the same curves. This value decreases 
till the point $\tilde y_i$, then it increases when the line
$\eta(t,x_1^i)$ departs from a neighborhood of the mountain pass point, to reach
a local minimizer. Remark that both figures represents just the first $1000$ steps in
the $t$-variable.

We now point out some cautions to be taken in the implementation of Algorithm \ref{alg:1b}.
After a certain number of steps, depending on the distance between $x_1$ and $x_2$, 
the points $x_1^i$ and $x_2^i$ may not be any more numerically distinct; 
moreover, taking into account the numerical errors in the integration method and a possible strong sharpness of the
graph of $f$, the lines of the steepest descent flow starting from $x_1^i$ and $x_2^i$ can become undistinguished.
Hence, each time we implement Algorithm \ref{alg:1b}
we a priori fix a maximum number of iterations, $N_{max}$, 
to avoid such numerical obstructions and to make
the distance $\mathrm{dist}(x_1^{N_{max}},x_2^{N_{max}}) \ll \epsilon$, where 
$\epsilon$ is a fixed small positive number.
Fixed $\epsilon,\gamma >0$ we propose the following algorithm that allows us 
to approach a locally optimal path joining the starting points $x_1$ and $x_2$.
\begin{algorithm}
\label{alg:1c}
\texttt{
\begin{itemize}
\item[Step 0] $N_{max}^0 = N_{max}$; \\
              \textsc{Algorithm 2}($x_1,x_2,N_{max}^0;\tilde y_{N_{max}^0}$)
\item[Step k] if $\| \nabla f(\tilde y_{N_{max}^{k-1}}) \| < \gamma $, STOP\\
              else $T_{def}:= \inf \left\{ t>0 : \mathrm{dist} 
              \left(\eta(t,x_1^{N_{max}^{k-1}}),\eta(t,x_2^{N_{max}^{k-1}})\right) \geq 
              \frac{\epsilon}{2^{k-1}}\right\}$;\\
              \hspace{1cm} $x_j := \eta(T_{def},x_j^{N_{max}^{k-1}})$, $j=1,2$;\\
              \hspace{1cm} $N_{max}^{k} := N_{max}^{k-1} -1$;\\
              \hspace{1cm} \textsc{Algorithm 2}($x_1,x_2,N_{max}^{k};\tilde y_{N_{max}^k}$)
\end{itemize}}
\end{algorithm}
It must be said that the choice of the second connected component where $x_2$ lies
is forced by the problem itself. Nevertheless, we can connect $x_1$ to any given connected 
component by juxtaposition of a finite number of locally optimal paths.

\subsection{Sublevels with non trivial fundamental group}
\label{subsec:not_sim}

We recall that a path $\gamma$ is a continuous function $\gamma : [0,1] \rightarrow X$.
Given a path $\gamma$, we define the path $-\gamma$ as 
$$
-\gamma(s) := \gamma(1-s), \quad \forall s \in [0,1]
$$
and for every pair of paths $\gamma_1,\gamma_2$ such that $\gamma_1(1)=\gamma_2(0)$ 
we define the path $\gamma = \gamma_1 \circ \gamma_2$ 
\emph{juxtaposition of $\gamma_1,\gamma_2$} as
$$
\gamma (s) = (\gamma_1 \circ \gamma_2)(s) := 
\left\{ \begin{array}{ll} \gamma_1(2s)   & s \in [0,\frac{1}{2}]\\
 													\gamma_2(2s-1) & s \in [\frac{1}{2},1]	
 			  \end{array}
\right.
$$
A path $\gamma : [0,1] \rightarrow X$ is a \emph{loop}  
if $\gamma(0) = \gamma(1)$. We say that a loop $\gamma$ is
\emph{contractible in $X$} if there exist ${\bar x} \in X$ and
a continuous function $h : [0,1] \times [0,1] \rightarrow X$
such that $h(0,s)=\gamma(s)$, $h(1,s)={\bar x}$, $\forall s \in [0,1]$, 
and $h(\lambda,0)=h(\lambda,1)$, $\forall \lambda \in (0,1)$.
We observe that, for every $\lambda \in [0,1]$, $h_\lambda(\cdot) := h(\lambda,\cdot)$ is 
a loop in $X$.

Let $Y$ be a subset of $X$, we say that $Y$ is \emph{simply connected} if every 
loop in $Y$ is contractible.

We consider a functional $f:X \rightarrow \RR$, 
satisfying the following hypotheses
\begin{itemize}
\item[{(h1)}] $f$ is bounded below, suppose $f(x) \geq 0, \quad \forall x \in X$;
\item[{(h2)}] there exists $c$ such that the sublevel $f^{c} \neq \emptyset$  is not simply connected;
\item[{(h3)}] $f$ verifies the (PS)$_{c_0}$ for every $0 \leq c_0 < c$.
\end{itemize}

In the sequel, we will work in a connected component of the sublevel $f^c$; we will still term it $f^c$.

Let $\bar x \in f^{c}$ be a strict local minimizer for $f$ and
$r>0$ such that the ball $B_r(\bar x)$ is contractible in $f^{c}$,
i.e. for every $x \in B_r(\bar x)$ the set 
$\{(1-\lambda)x + \lambda\bar x : \lambda \in [0,1]\} \subset f^{c}$.
We consider the sublevel $f^{\bar c + \epsilon}$, where 
$\bar c = f(\bar x)$ and $\epsilon > 0$, such that $F^{\bar c + \epsilon}_{\bar x} \subset B_r(\bar x)$,
where $F^{\bar c + \epsilon}_{\bar x}$ is the connected component of
$f^{\bar c + \epsilon}$ containing ${\bar x}$. 
We observe that $F^{\bar c + \epsilon}_{\bar x}$ is contractible in $f^{c}$.

For every $x_1 \in f^{c}$ such that $\omega_{x_1} = \bar x$ 
we define the instant 
$$
T_{x_1} := \inf \{ t\geq 0 : \eta(t,x_1) \in F^{\bar c + \epsilon}_{\bar x} \} > 0.
$$
We define the following paths
\begin{equation}
\label{desc_path_12}
\gamma^1_{x_1}(\lambda) := \eta(T_{x_1}\lambda,x_1), \quad
\gamma^2_{x_1}(\lambda) := (1-\lambda)\eta(T_{x_1},x_1) + \lambda{\bar x}, 
\end{equation}

\begin{definition}
\label{def:desc_path}
We term \emph{descending path associated to $x_1$} the path 
$$
\alpha_{x_1} := (\gamma^1_{x_1} \circ \gamma^2_{x_1}) 
$$
where $\gamma^1_{x_1}, \gamma^2_{x_1}$ are defined in~(\ref{desc_path_12}).
\end{definition}

\begin{remark}
\label{rem:cont}
We observe that, defining $\bar X :=\{x \in X : \omega_x =\bar x \}$
and $\Gamma_{\bar X, \bar x} := \{ \gamma \in \Gamma_{x, \bar x} : x \in \bar X\}$,
the maps
$$
\bar X \rightarrow \RR_+, \, x_1 \mapsto T_{x_1}  
\quad \mbox{ and }\quad
\bar X \rightarrow \Gamma_{\bar X, \bar x}, \,  x_1 \mapsto \alpha_{x_1} 
$$
are continuous on the set $\bar X$.
\end{remark}

\begin{definition}
\label{def:desc_loop}
Let $\gamma$ be a path in $f^{c}$ such that 
$\omega_{\gamma(0)} = \omega_{\gamma(1)} = {\bar x}$; 
we define the loop $\gamma_\eta$, the \emph{descending loop associated to $\gamma$}, as
\begin{equation}
\label{desc_loop}
\gamma_\eta := [(-\alpha_{\gamma(0)}) \circ \gamma] \circ \alpha_{\gamma(1)}.
\end{equation}
\end{definition}

\begin{definition}
\label{def:eta_contr}
We say that a path $\gamma$ in $f^{c}$
is \emph{$\eta$-contractible in $f^{c}$} if 
\begin{itemize}
\item[{(a)}] $\omega_{\gamma(0)} = \omega_{\gamma(1)} = {\bar x}$;
\item[{(b)}] the loop $\gamma_\eta$ is contractible in $f^{c}$.
\end{itemize}
\end{definition}

\begin{lemma}
\label{lem:eta_contr}
Let $\bar x \in f^{c}$ be a strict local minimizer for the functional $f$
and $\gamma$ a path in $f^{c}$. If $\omega_{\gamma(s)}={\bar x}$, $\forall s \in [0,1]$,
then $\gamma$ is $\eta$-contractible in $f^{c}$.
\end{lemma}

\begin{proof}
To prove that the path $\gamma_\eta$ defined in~(\ref{desc_loop})
is contractible, we define a continuous function 
$h: [0,1] \times [0,1] \rightarrow f^{c}$ such that
$h(0,\cdot)=\gamma_\eta(\cdot)$ and $h(1,\cdot)=\bar x$.
For every  $\lambda \in \left[\frac{1}{2},1\right]$
we define the path
$$
\beta_\lambda(s) := \alpha_{\gamma(1)}(2(1-\lambda)s), \quad s \in \left[0,1\right],
$$
hence, the following function $h$ satisfies our requirements 
\begin{equation*}
h_\lambda := \left\{ 
\begin{array}{ll}\displaystyle
[(-\alpha_{\gamma(2\lambda)}) \circ \gamma] \circ \alpha_{\gamma(1)}, 
& \lambda \in \left[0,\frac{1}{2}\right] \\
(-\beta_\lambda) \circ \beta_\lambda,                                        
& \lambda \in \left[\frac{1}{2},1\right]
\end{array} \right.
\end{equation*}
\end{proof}
 
\begin{theorem}
\label{thm:mp2}
Let $f$ be a $C^2$ functional on a Hilbert space $X$.
Suppose that $f$ satisfies conditions \textup{{(h1)}, {(h2)}, {(h3)}}.
Then one of the following situations occurs:
\begin{itemize}
\item[(i)] $f^{c}$ has a continuum of global minimizers;
\item[(ii)] there exists ${\hat x} \in Crit(f) \cap f^{c}$ that is not a strict minimizer.
\end{itemize}
\end{theorem}

\begin{proof}
When the set $Crit(f) \cap f^{c}$ is not a continuum of global minimizers, 
hypotheses {(h1)} and {(h3)} imply the existence of at least one strict minimizer 
$\bar x \in f^{c}$ for the functional $f$. 
Let $\gamma \subset f^{c}$ be a loop  such that 
$\gamma(0)=\gamma(1)=\bar x$ and suppose that $\gamma$ is not contractible in $f^{c}$.
Hence, since $\gamma_\eta(s)=\gamma(s)$ for every $s \in [0,1]$, 
$\gamma$ is not $\eta$-contractible in $f^{c}$.
Using Lemma~\ref{lem:eta_contr}, we can conclude that
there exists $x_\gamma \in \gamma([0,1])$ such that $\omega_{x_\gamma} \neq \bar x$.
The critical point $\omega_{x_\gamma}$ could be a strict minimizer or not. In the first case, 
we use Algorithm~\ref{alg:1b} to prove the existence of a critical point for $f$ that is not
a strict local minimizer; otherwise $\omega_{x_\gamma}$ is the searched critical point for $f$.

We now show an algorithm that determine a point $x_\gamma \in \gamma([0,1])$
such that $\omega_{x_\gamma} \neq \bar x$. 
Let $\gamma \subset f^{c}$ be a not contractible loop such that 
$\gamma(0)=\gamma(1)=\bar x$. We have already remarked that $\gamma$ is
$\eta$-contractible. 

\begin{algorithm}
\label{alg:2}
\texttt{
\begin{itemize}
\item[{ Step 0.}] $\gamma_0 := \gamma$, $x_0 = \gamma_0\left( \frac{1}{2} \right)$\\
                  if $\omega_{x_0} \neq \bar x$, STOP\\
                  if $\gamma_0\left(\left[ 0,\frac{1}{2} \right]\right)$ is $\eta$-contractible,
                  $\gamma_1 := \gamma_0\left(\left[\frac{1}{2},1\right]\right)$\\
                  else $\gamma_1 := \gamma_0\left(\left[0,\frac{1}{2}\right]\right)$\\
\item[{ Step i.}] $x_i = \gamma_i\left( \frac{1}{2} \right)$\\
                  if $\omega_{x_i} \neq \bar x$,  STOP\\
                  if $\gamma_i\left(\left[0,\frac{1}{2}\right]\right)$ is $\eta$-contractible,
                  $\gamma_{i+1} := \gamma_i\left(\left[\frac{1}{2},1\right]\right)$\\
                  else $\gamma_{i+1} := \gamma_i\left(\left[0,\frac{1}{2}\right]\right)$.
\end{itemize}}
\end{algorithm}
Using Algorithm~\ref{alg:2} we can directly find a point whose $\omega$-limit 
is not $\bar x$, in this case the $i$-loop is stopped. 
Otherwise, we define a sequence of paths $(\gamma_i)_i$, that are not $\eta$-contractible. In this case, let $x_i^0=\gamma_i(0)$ and $x_i^1=\gamma_i(0)$
be the initial and end points of these paths, then the sequences
$(x_i^0)_i, (x_i^1)_i \subset \gamma([0,1])$ and 
$$
\lim_{i \rightarrow +\infty} x_i^0 = \lim_{i \rightarrow +\infty} x_i^1 = \bar{\bar x}.
$$
Arguing as in the proof of the convergence of 
Algorithm~\ref{alg:1a}, we deduce that necessarily $\omega_{\bar{\bar x}} \neq \bar x$;
if $\omega_{\bar{\bar x}}$ is a strict minimizer, we use Algorithm~\ref{alg:1b} to deduce the 
existence of a critical point that is not a strict local minimizer.
\end{proof}

\section{The Mountain Pass Theorem}
\label{sec:2}

Our goal now is to prove a version of the Mountain Pass Theorem 
(for a detailed theory on this subject we refer to~\cite{AR,PS1}).
A Mountain Pass Theorem concerns itself with proving the existence of 
critical points which are not strict local minimizers for the functional $f$; 
using Corollary~\ref{cor:thm:alg1} and Proposition~\ref{prop cal F}, we are now able to prove the following

\begin{theorem}[Mountain Pass Theorem]
\label{thm:mp1}
Let $f$ be a $C^2$ functional on  a Hilbert space $X$.
Let $x_1,x_2 \in X$, let $\Gamma_{x_1,x_2}$ be the set of paths defined in~\textup{(\ref{Gamma_set})} and
$c_0$ the level
\begin{equation}
\label{bar c}
{c}_0:= \inf_{\gamma \in \Gamma_{x_1,x_2}} \sup_{s \in [0,1]}f(\gamma(s)).  
\end{equation}
such that 
\begin{equation}
\label{c>max}
{c}_0 > \max \{f(x_1),f(x_2)\}.
\end{equation}
If the functional $f$ satisfies the Palais-Smale condition at level ${c_0}$
then there exists a critical point for the functional $f$ at level $c_0$,
that is not a local minimizer. 
\end{theorem} 

\begin{proof}
The definition~(\ref{bar c}) of the level $c_0$ and condition~(\ref{c>max}), 
imply, first, that the sublevel $f^{c_0}$ is disconnected, 
second that for every $k \in \NN$,
there exist $\gamma_k \in \Gamma_{x_1,x_2}$ such that 
$$
\sup_s f(\gamma_k(s)) \leq c_0 + \frac{1}{k}.
$$
From Corollary~\ref{cor:thm:alg1}, we deduce the existence of an element
$\bar x_k \in \gamma_k([0,1]) \cap \partial \mathcal{F}^{c_0}_1$ such that
$$
c_0 \leq f(\omega_{{\bar x}_k}) \leq c_0 + \frac{1}{k}.
$$
Let $(x_{1,k}^n)_n \subset \mathcal{F}^{c_0}_1$ converging to $\bar x_k$;
following the proof of Corollary~\ref{cor2:thm:alg1}, we deduce the 
existence of an index $n_k \in \NN$ and an instant $\tilde T_{n_k}$ such that, 
$\tilde y_{n_k} := \eta (\tilde T_{n_k},x_1^{n_k})$, where
$x_1^{n_k} \in (x_{1,k}^n)_n$, verifies
$$
\| \nabla f(\tilde y_{n_k}) \| < \frac{1}{k} \mbox{ and }
c_0 \leq f(\tilde y_{n_k}) \leq c_0 +\frac{1}{k}.
$$
Since $f$ satisfies the Palais-Smale condition at level ${c_0}$, we deduce that 
there exists a critical point, $\tilde y$, for $f$ at level $c_0$, 
limit of the sequence $(\tilde y_{n_k})_{k}$ as $k$ tends to $+\infty$.

We conclude the proof showing that there exists a sequence $(y^*_{n_k})_k$, converging to 
a critical point $y^*$ that is not a minimizer and such that 
$f(y^*_{n_k}) \rightarrow f(y^*)=c_0$, with $f(y^*_{n_k}) = c_0 -\frac{1}{k} <c_0$.
Let 
$$
T^*_{n_k} := \inf \left\{ t \geq 0 : f(\eta(t,x^{n_k}_1)) \leq c_0 -\frac{1}{k} \right\},
$$
hence, defining $y^*_{n_k} := \eta(T^*_{n_k},x^{n_k}_1)$, we have $f(y^*_{n_k}) = c_0 -\frac{1}{k}$.
If we can prove that 
\begin{equation*}
\lim_{k \rightarrow +\infty} \nabla f(y^*_{n_k}) = 0,
\end{equation*}
using the Palais-Smale condition at level $c_0$ for $f$, we conclude the proof.
By the sake of contradiction suppose that there exists $\epsilon_0 >0$
such that $\| \nabla f(y^*_{n_k}) \| \geq \epsilon_0$, for every $k \in \NN$.
Hence, since the Palais-Smale condition at level $c_0$ holds,
there exists $\delta>0$ such that $|T^*_{n_k} - \tilde T_{n_k}| > \delta$ 
and $\| \nabla f(\eta(T^*_{n_k}-t,x^{n_k}_1))\| \geq \epsilon' > 0$, $\forall t \in [0,\delta)$.
Then, using Lemma~\ref{mis_tempi}
\begin{multline*}
0< \epsilon' \leq \inf_{t \in [0,\delta)} \| \nabla f(\eta(T^*_{n_k}-t,x^{n_k}_1))\| 
          \leq \sqrt{\frac{2(f(\eta(T^*_{n_k}-\delta,x^{n_k}_1))-f(y^*_{n_k}))}{\delta^2}} \\
          \leq \sqrt{\frac{2\left(f(\tilde y_{n_k})-c_0 + \frac{1}{k}\right)}{\delta^2}} \leq 
          \frac{1}{\sqrt{k\delta}},          
\end{multline*}
and the right-hand side tends to $0$ as $k$ tends to $+\infty$.
\end{proof}

\bibliographystyle{amsplain}

\begin{thebibliography}{10}
%
\bibitem{AR}
{A. Ambrosetti and P.H. Rabinowitz}, 
{\it Dual variational methods in critical point theory and applications},
J. Funct. Anal. {\bf 14} (1973), 349--381. 
%
\bibitem{PhD_vi}
{V. Barutello},
{\it On the $n$-body problem},
Ph.D thesis, Universit\`a di Milano-Bicocca, Milano 2004.
%
\bibitem{CMc}
{Y.S. Choi and P.J. McKenna}, 
{\it  A mountain pass method for the numerical solution of semilinear elliptic problems}, 
Nonlinear Anal. {\bf 20} 4 (1993), 417--37.
%
\bibitem{CMcR}
{Y.S. Choi, P.J. McKenna and M. Romano}, 
{\it  A mountain pass method for the numerical solution of semilinear wave equations}, 
Numer. Math. {\bf 64} (1993), 487--459.
%
\bibitem{PS1}
{P. Pucci and J. Serrin},
{\it The structure of the critical set in the Mountain Pass Theorem},
Trans. Amer. Math. Soc. {\bf 299} 1 (1987), 1115--132.
%
\end{thebibliography}

\end{document}